\documentclass{article}

\newtheorem{fed}{\textbf{Definition}}[section]
\newtheorem{thm}[fed]{\textbf{Theorem}}

\usepackage{amssymb,bbm,graphicx,epsfig,psfrag,epic,eepic,latexsym}
\usepackage{amsmath}
\usepackage{mathrsfs}
\usepackage[dvips]{color}

\begin{document}
\title{Dihedral homology and the moon}
\author{Urs Frauenfelder\footnote{Department of Mathematics and
Research Institute of Mathematics, Seoul National University}}
\maketitle

\begin{abstract}
The equivariant homology of a twisted action of $O(2)$ on the free
loop space of spheres was computed by Lodder \cite{lodder} via
dihedral homology. In this note we explain how this twisted action
appears in the restricted three body problem after Moser
regularization.
\end{abstract}

\section{Introduction}

If $N$ is a manifold, the standard action of the group $O(2)$ on the
circle $S^1$ gives rise to an $O(2)$-action on the free loop space
$\mathscr{L}_N=C^\infty(S^1,N)$ by
$$g_*v(t)=v(gt), \quad v \in \mathscr{L}_N,\,\,g \in O(2),\,\,t \in
S^1.$$ This $O(2)$-action plays a major role in the search of closed
geodesics. Indeed, if $N=S^2$ is the two dimensional sphere, then
this $O(2)$-action is the basic reason for the existence of the
three simple closed Lusternik-Schnirelman geodesics
\cite{ballmann,lusternik-schnirelmann}. How dramatic the situation
changes without the $O(2)$-action is impressively illustrated by the
Katok examples \cite{ziller}. Indeed, for nonreversible Finsler
metrics the geodesic equation is only invariant under the action of
$SO(2)$, and as the Katok examples illustrate, there are examples of
nonreversible Finsler metrics on $S^2$ with only two closed
geodesics.

However, if we denote by $\mathscr{L}$ the free loop space of the
$n$-dimensional sphere $S^n=\{p \in \mathbb{R}^{n+1}: ||p||=1\}$ we
can consider the following twisted $O(2)$-action on $\mathscr{L}$.
Let
\begin{equation}\label{refl}
\rho \colon S^n \to S^n, \quad (p_1,p_2 \ldots, p_{n+1}) \mapsto
(-p_1,p_2 \ldots, p_{n+1})
\end{equation}
be the reflection at the equator. For $g \in O(2)$ set
$$\iota(g)=\frac{1-\det g}{2},$$
i.e. $\iota(g)=0$ if $g \in SO(2)$ and $\iota(g)=1$ if $g \in O(2)
\setminus SO(2)$. The twisted $O(2)$-action on $\mathscr{L}$ is
given for $g \in O(2)$ and $v \in \mathscr{L}$ by
$$g_*v(t)=\rho^{\iota(g)}v(gt), \quad t \in S^1.$$
Since $\rho$ is an involution this defines indeed an action.

For the twisted-$O(2)$ action the equivariant homology of
$\mathscr{L}$ was completely computed by Lodder \cite{lodder} after
relating it to dihedral homology \cite{fiedorowicz-loday, loday}.
Moreover, Lodder constructed in this paper a combinatorial model for
the homotopy quotient which also appears in the work of
B\"odigheimer-Madsen \cite{bodigheimer-madsen}. In contrast to the
twisted case to the authors knowledge we only partially know the
equivariant homology of the standard $O(2)$-action on $\mathscr{L}$
due to work of Hingston \cite{hingston1}.

From the point of view of closed geodesics it might seem a bit
strange to look at the twisted $O(2)$-action on $\mathscr{L}$.
Indeed, to take advantage of this action one has to restrict one's
attention to metrics which are invariant under the reflection
$\rho$. However, the purpose of this note is to explain that the
twisted $O(2)$-action actually naturally arises in the study of a
Hamiltonian system whose puzzling dynamics historically was one of
the main motivations to investigate closed geodesics - namely the
restricted problem of three bodies \cite{barrow-green}. Let us see
what Poincar\'e himself wrote about this topic \cite{poincare}.
\begin{quote}
Dans mes M\'ethodes nouvelles de la M\'ecanique c\'eleste j'ai
\'etudi\'e les particularit\'es des solutions du probl\`eme des
trois corps et en particulier des solutions p\'eriodiques et
asymptotiques. Il suffit de se reporter \`a ce que j'ai \'ecrit \`a
ce sujet pour comprendre l'extr\^eme complexit\'e de ce probl\`eme;
\`a c\^ot\'e de la difficult\'e principale, de celle qui tient au
fond m\^eme des choses, il y a une foule de difficult\'es
secondaires qui viennent compliquer encore la t\^ache du chercheur.
Il y aurait donc int\'er\^et \`a \'etudier d'abord un probl\`eme
o\`u on rencontrerait cette difficult\'e principale, mais o\`u on
serait affranchi de toutes les difficult\'es secondaires. Ce
probl\`eme est tout trouv\'e, c'est celui des lignes g\'eod\'esiques
d'une surface.
\end{quote}
Contrary to Poincar\'e's expections this note indicates that at
least some aspects of the second problem are even harder than they
were in the original one.

To interpret the flow of the restricted three body problem as a
Hamiltonian flow on the cotangent bundle of a sphere one first has
to regularize collisions. A beautiful way how two body collisions
can be regularized was described by Moser \cite{belbruno, moser}. If
one interprets $S^n=\mathbb{R}^n \cup \{\infty\}$, then the point at
infinity corresponds to collisions.

Trajectories of the restricted three body problem are usually not
invariant under reversal of time. Therefore the standard
$O(2)$-action on $\mathscr{L}$ is of no use for the restricted three
body problem. On the other hand at least since the fundamental work
of Birkhoff \cite{birkhoff} an important aspect in the study of the
dynamics of the restricted three body problem is the fact that its
Hamiltonian is invariant under a different antisymplectic
involution. This involution extends to the regularization and hence
gives rise to an $O(2)$-action on periodic trajectories of the
restricted three body problem. In this note we explain the following
observation.
\\ \\
\textbf{Observation: } \emph{After Moser regularization the
$O(2)$-action on periodic trajectories of the circular restricted
three body problem coincides with the twisted $O(2)$-action on
$\mathscr{L}$ described above.}
\\ \\
In \cite{albers-frauenfelder-koert-paternain} it was shown that in
the planar case, i.e. $n=2$, below the first critical value, the
bounded components of the energy hypersurfaces of the regularized
restricted three body problem in $T^*S^2$ are fiberwise starshaped.
This implies that their symplectic homology is well defined and
coincides with the loop space homology of $S^2$
\cite{abbondandolo-schwarz, salamon-weber, viterbo}. In particular,
if one looks at the $O(2)$-equivariant symplectic homology defined
as in \cite{bourgeois-oancea}, one gets the following Corollary if
one combines Lodder's computation \cite{lodder} with the
Observation.
\\ \\
\textbf{Corollary: } \emph{Below the first critical value the
$O(2)$-equivariant symplectic homology for each of the bounded
components of the regularized circular planar restricted 3-body
problem is given by the following formula}
$$SH_*^{O(2)}=H_*(BO(2)) \oplus \bigoplus_{m=1}^\infty
H_*(D_m;\mathbb{Z}[1]^{ \otimes m}).$$ \textbf{Remark: } It is
conceivable that the methods in
\cite{albers-frauenfelder-koert-paternain} can also be used to show
that in the spatial case of the restricted three body problem, i.e.
$n=3$, below the first critical value the bounded components of the
regularized energy hypersurface are still fiberwise starshaped in
$T^*S^3$. Hence the $O(2)$ equivariant symplectic homology would
still be well defined and could be computed as in the planar case to
be $SH_*^{O(2)}=H_*(BO(2)) \oplus \bigoplus_{m=1}^\infty
H_*(D_m;\mathbb{Z}[2]^{\otimes m})$.
\\ \\
\textbf{Remark: } It is an open question if the bounded components
of the energy hypersurfaces for the regularized restricted three
body problem are also fiberwise convex. If this were true, then
after regularization trajectories could be interpreted as Finsler
geodesics for a (nonreversible) Finsler metric. In this case
symplectic homology in the Corollary above could be replaced by the
Morse homology of a Finsler metric. It was shown in
\cite{cieliebak-frauenfelder-koert} that in the special case of the
rotating Kepler problem, i.e. the case were the mass of the moon
becomes zero, fiberwise convexity actually holds.
\\ \\
This note is organized as follows. In Section~\ref{dh} we recall
Lodder's computations of the equivariant homology of the loop space
of a sphere for the twisted $O(2)$-action. In Section~\ref{mr} we
explain how the Kepler flow can be interpreted as the geodesic flow
of a round sphere after Moser regularization. In Section~\ref{rp} we
explain the restricted three body problem and its Moser
regularization. In Section~\ref{bi} we explain how Birkhoff's
involution extends to the Moser regularization and prove the
Observation and its Corollary.
\\ \\
\emph{Acknowledgements: } The author would like to thank Nancy
Hingston for drawing his attention to Lodder's paper \cite{lodder}.
This research was partially supported by the Basic Research fund
2010-0007669 funded by the Korean government.

\section{Dihedral homology}\label{dh}

To the authors knowledge historically the first complete computation
of the $SO(2)$-equivariant homology of free loop spaces of spheres
appeared in the paper by Carlsson and Cohen \cite{carlsson-cohen}
and relied on deep results from cyclic homology. Before the paper by
Carlsson and Cohen partial computations were achieved by Hingston
\cite{hingston1} using Morse theory. An interesting Morse theoretic
interpretation of the result by Carlsson and Cohen as well as an
alternative derivation of their results can be found in
\cite{hingston2}.

Since the scheme of Lodders computation of the twisted
$O(2)$-equivariant homology of spheres follows the scheme of
Carlsson and Cohen we first recall their computation. It was proved
by Burghelea-Fiedorowicz, Goodwillie and Jones that if $N$ is a
connected manifold and $\mathscr{L}_N$ the free loop space of $N$,
then
\begin{equation}\label{goodwil}
H_*^{SO(2)}(\mathscr{L}_N)=HC_*(C_*(\Omega N)). \end{equation} Here
$\Omega N$ is the Moore loop space of $N$ which is a group like
topological monoid, $C_*(\Omega N)$ is the differential graded
algebra of singular chains of $\Omega N$, and $HC_*$ is the cyclic
hyperhomology of this differential graded algebra. Strictly speaking
the paper of Jones \cite{jones} relates the equivariant loop space
homology to a variant of cyclic hyperhomology applied to the
differential graded algebras of cochains on $N$. How the two
approaches are related is explained in \cite{jones-mccleary}.

If $n \geq 2$, then the Pontryagin ring of $\Omega S^n$ is given by
$$H_*(\Omega S^n)=\mathbb{Z}[x], \quad |x|=n-1.$$
If we think of $H_*(\Omega S^n)$ as a differential graded algebra
with trivial differential, then we can construct a quasiisomorphism
between $H_*(\Omega S^n)$ and $C_*(\Omega S^n)$ by mapping $x$ to a
representative of it. Since cyclic hyperhomology is invariant under
quasiisomorphism we conclude for $\mathscr{L}=\mathscr{L}_{S^n}$
that
$$H_*^{SO(2)}(\mathscr{L})=HC_*(\mathbb{Z}[x]), \quad |x|=n-1.$$
The latter was computed by Loday and Quillen \cite{loday,
loday-quillen} to be
$$H_*^{SO(2)}(\mathscr{L})=H_*(BSO(2))\oplus \bigoplus_{m=1}^\infty
H_*(\mathbb{Z}_m;\mathbb{Z}[n-1]^{\otimes m})$$ where $BSO(2) \cong
CP^\infty$ is the classifying space of the group $SO(2)$.

Thanks to the theory of crossed simplicial groups due to Fiedorowicz
and Loday \cite{fiedorowicz-loday}, there are analogons of
(\ref{goodwil}) for more general actions then circle actions. In
particular,
\begin{equation}\label{dunn}
H^{O(2)}_*(\mathscr{L})=HD_*(C_*(\Omega S^n))
\end{equation}
where on the righthand side we have dihedral hyperhomology
\cite{dunn}. This is actually true for both the standard and the
twisted $O(2)$-action on $\mathscr{L}$, where we have to note that
the action of the dihedral group on the Hochschild complex of
$C_*(\Omega S^n)$ differs in both cases. For a notation which keeps
track of the actions we refer again to Dunn's paper \cite{dunn}.

However, just for the twisted action we can identify the righthand
side of (\ref{dunn}) with the dihedral homology of the Pontrjagin
ring of $\Omega S^n$. First pick a basepoint on $S^n$ and identify
$S^n$ with the reduced suspension $S^n=\Sigma S^{n-1}=S^{n-1} \wedge
S^1$. Consider the map $\alpha \colon S^{n-1} \to \Omega S^n$ given
for $p \in S^{n-1}$ by
$$\alpha(p)(t)=(p,t) \in S^{n-1} \wedge S^1, \quad t \in [0,1].$$
Choose a cycle $\sigma \in C_{n-1} S^{n-1}$ representing a homology
class $[\sigma]$ which generates $H_{n-1}(S^{n-1})\cong \mathbb{Z}$.
Then $\alpha_*[\sigma] \in H_{n-1}(\Omega S^n)$ is a generator of
the Pontrjagin ring. The cycle $\alpha_*\sigma \in C_{n-1}(\Omega
S^n)$ is invariant, if one simultaneously switches the suspension
coordinate as well as the direction of a loop in $\Omega S^n$. Hence
for the twisted $O(2)$-action on $\mathscr{L}$ one gets from
(\ref{dunn})
$$H^{O(2)}_*(\mathscr{L})=HD_*(\mathbb{Z}[x]), \quad |x|=n-1.$$
The dihedral homology of $\mathbb{Z}[x]$ was computed by Lodder in
\cite{lodder} to be
\begin{equation}\label{lod}
H^{O(2)}_*(\mathscr{L})=H_*(BO(2))\oplus \bigoplus_{m=1}^\infty
H_*(D_m;\mathbb{Z}[n-1]^{\otimes m})
\end{equation}
where $D_m$ is the dihedral group of order $2m$.
\\ \\
\textbf{Remark: } The computations of Carlsson and Cohen
\cite{carlsson-cohen} as well as of Lodder \cite{lodder} which we
recalled for spheres actually apply more generally to suspensions.

\section{Moser regularization}\label{mr}

The Hamiltonian for the Kepler problem $H \colon T^*(\mathbb{R}^n
\setminus \{0\})=(\mathbb{R}^n\setminus \{0\}) \times \mathbb{R}^n
\to \mathbb{R}$ is given by
$$H(q,p)=\frac{1}{2}|p|^2-\frac{1}{|q|}.$$
For negative energy values trajectories of the Hamiltonian flow are
either collision orbits are project to ellipses in position space.
In particular, for a dense set of trajectories the flow is periodic
and the only nonperiodic orbits are collision orbits. Two body
collisions can always be regularized and in \cite{moser} Moser found
a way to embed a reparametrization of the Kepler flow into the
geodesic flow of the round metric of $S^n$. In particular, after
regularization the flow becomes completely periodic.

Moser's method is most nicely illustrated for the energy value
$c=-\frac{1}{2}$. Consider the Hamiltonian $K \colon
(\mathbb{R}^n\setminus \{0\}) \times \mathbb{R}^n \to \mathbb{R}$
defined by
$$K(q,p)=|q|\bigg(H(q,p)+\frac{1}{2}\bigg)+1=\frac{1}{2}(|p|^2+1)|q|.$$
Note that
$$\Sigma:=H^{-1}\bigg(-\frac{1}{2}\bigg)=K^{-1}(1).$$
For a point $(q,p) \in \Sigma$ the Hamiltonian vector fields satisfy
$$X_K(q,p)=|q|X_H(q,p).$$
Since $|q|$ never vanishes on the energy hypersurface $\Sigma$, the
flow of $X_K$ restricted to $\Sigma$ is just a reparametrization of
the flow of $X_H$ on $\Sigma$. Now switch in your mind the roles of
$q$ and $p$ and think of $p$ as the base coordinate and $q$ as the
fiber coordinate. Moreover, think of $\mathbb{R}^n$ as the chart of
$S^n$ obtained via stereographic projection. Then $K(q,p)$ is
precisely the length of the cotangent vector $q \in T^* S^n$ for the
round metric on $S^n$. Note further that the switch of base and
fiber coordinates $(q,p) \mapsto (-p,q)$ is a symplectomorphism from
$T^*\mathbb{R}^n$ to $T^*\mathbb{R}^n$. Hence the Hamiltonian flow
of $K$ embeds into the geodesic flow of the round metric on $S^n$
and therefore after reparametrization the Kepler flow as well.

If one thinks of the sphere as $S^n=\mathbb{R}^n \cup \{\infty\}$,
then the point at infinity corresponds to collisions. Here it is
useful to remember that we switched the roles of position and
momentum coordinates. Indeed, at collisions the original momentum
explodes while the original position coordinate remains bounded.

If the energy value $c$ is negative but not necessarily equal to
$-\frac{1}{2}$ we still get the geodesic flow on $S^n$ after
conjugating with an additional diffeomorphism in the chart obtained
via stereographic projection. Indeed, for negative $c$ set
$$K_c(q,p)=|q|\big(H(q,p)-c\big)+1=\frac{1}{2}(|p|^2+|c|)|q|.$$
Consider the symplectomorphism $\phi \colon T^* \mathbb{R}^n \to T^*
\mathbb{R}^n$ given by
$$\phi(p,q)=\bigg(\sqrt{|c|} p,\frac{q}{\sqrt{|c|}}\bigg), \quad (p,q)
\in T^* \mathbb{R}^n.$$ Then we get
$$\phi^* K_c(p,q)=\frac{\sqrt{|c|}}{2}(|p|^2+1)|q|$$
which up to a constant conformal factor is again the length of the
cotangent vector $q$ for the standard round metric on $S^n$ in the
chart obtained by stereographic projection.
\section{The restricted three body problem}\label{rp}

In the restricted three body problem one considers two massive
bodies, the primaries, and a massless body which is attracted by the
two primaries according to Newton's law of gravitation. We refer to
the primaries as the earth and moon and to the massless body as the
satellite. Since the satellite is assumed to be massless it does not
affect the motion of the two primaries. Hence the motion of the two
primaries is governed by Kepler's laws. In the \emph{circular}
restricted three body problem one assumes in addition that the two
primaries move along circles around their common center of mass. One
further distinguishes between the \emph{planar} case, where the
satellite is supposed to move in the eccliptic, i.e. the plane
spanned by earth and moon, and the \emph{spatial} case where the
satellite is supposed to move in three dimensional space.
Mathematically the problem makes sense in any dimension greater or
equal to two and we therefore describe the Hamiltonian for the
movement of the satellite in a space of arbitrary dimension $n \geq
2$.

After scaling the total mass to be one and excluding the case that
the mass of the earth vanishes there is a $\mu \in [0,1)$ such that
the mass of the moon equals $\mu$ and the mass of the earth equals
$1-\mu$. Let $\{e_1, \cdots, e_n\}$ be the standard base of
$\mathbb{R}^n$. Scaling the distance of the earth and the moon to
one, translating the center of mass to the origin and applying a
suitable orthogonal transformation in the inertial system the
position of earth $E_i(t)$ and moon $M_i(t)$ at time $t$ are given
by
$$E_i(t)=\mu \cos(t) e_1+\mu \sin(t) e_2, \quad
M_i(t)=-(1-\mu) \cos(t) e_1-(1-\mu)\sin(t) e_2.$$ The Hamiltonian
for the satellite in the inertial system
$$H^i_t \colon T^*(\mathbb{R}^n \setminus \{E_i(t),M_i(t)\}) \to \mathbb{R}$$
is given by kinetic and potential energy
$$H^i_t(q,p)=\frac{1}{2}|p|^2-\frac{1-\mu}{|q-E_i(t)|}-
\frac{\mu}{|q-M_i(t)|}.$$ Note that since the earth and moon are
moving in the inertial system this Hamiltonian is not autonomous,
with, alas, even a time dependent domain of definition. In
particular, it is not preserved along trajectories of the satellite.
To improve this unpleasant situation we transform the system form
inertial to rotating coordinates. In the rotating coordinate system
the positions of the earth and moon are fixed
$$E=\mu e_1, \quad M=-(1-\mu) e_1.$$
The transition from the inertial to the rotating coordinate system
involves a time dependent transformation, namely we have to rotate
the $(e_1,e_2)$-plane in $\mathbb{R}^n$. The infinitesimal generator
of this rotation is given by angular momentum
$$L(q,p)=p_1 q_2-p_2 q_1.$$
Hence the Hamiltonian $H \colon T^*(\mathbb{R}^n \setminus \{E,M\})
\to \mathbb{R}$ in the rotating coordinate system becomes
\cite{abraham-marsden}
\begin{equation}\label{hami}
H(q,p)=\frac{1}{2}|p|^2-\frac{1-\mu}{|q-E|}-
\frac{\mu}{|q-M|}+L(q,p).
\end{equation}
Note that in the rotating coordinate system the Hamiltonian is
autonomous and therefore preserved under its Hamiltonian flow. This
would in general not be true for the elliptic restricted three body
problem where the primaries are allowed to move on ellipses. For
that reason it is important to restrict to the circular case. Since
this observation goes back to Jacobi the integral of motion $-2H$ is
usually referred to as the Jacobi integral.

Introducing the \emph{effective potential} $U \colon \mathbb{R}^n
\setminus \{E,M\} \to \mathbb{R}$ given by
$$U(q)=-\frac{1-\mu}{|q-E|}-\frac{\mu}{|q-M|}-\frac{1}{2}(q_1^2+q_2^2)$$
the Hamiltonian $H$ can be rewritten as
$$H(q,p)=\frac{1}{2}\bigg((p_1+q_2)^2+(p_2-q_1)^2+\sum_{i=3}^n
p_i^2\bigg)+U(q).$$ Note that the footpoint projection $\pi:T^*
(\mathbb{R}^n \setminus \{E,M\}) \to \mathbb{R}^n \setminus \{E,M\}$
restricts to a bijection between critical points
$$\Pi=\pi|_{\mathrm{crit}(H)} \colon \mathrm{crit}(H) \to \mathrm{crit}(U).$$
If $\mu \neq 0$, i.e.\,the mass of the moon does no vanish, there
are five critical points of $U$ for every massratio, usually
referred to as \emph{Lagrange points}. The first three Lagrange
points $L_1$, $L_2$ and $L_3$ are saddle points of $U$ and are
collinear with the earth and moon, $L_1$ lies between $E$ and $M$,
$L_2$ lies to the right of $E$, and $L_3$ lies to the left of $M$.
The Lagrange points $L_4$ and $L_5$ are maxima of $U$. They lie in
the $(e_1,e_2)$-plane and together with each of the primaries span
an equilateral triangle. In the limiting case where $\mu$ is zero,
i.e. the rotating Kepler problem, the critical set of $U$ consists
of the circle of radius one around the origin in the
$(e_1,e_2)$-plane.

For an energy value $c \in \mathbb{R}$ abbreviate by
$\Sigma_c=H^{-1}(c)$ the energy hypersurface of $H$. The
\emph{Hill's region} is defined to be
$$\mathcal{K}_c=\pi(\Sigma_c)=\big\{q \in \mathbb{R}^n \setminus \{E,M\}:
U(q) \leq c\big\}.$$ When $\mu \in (0,1)$ it can be shown that $L_1$
has the smallest action of all five critical points of $U$. If
$c<U(L_1)=H(\Pi^{-1}(L_1))$ the Hill's region has three connected
components. Two of these connected components are bounded. One
bounded component contains the earth in its closure and one bounded
component contains the moon in its closure. The third region is
unbounded. We denote by $\mathcal{K}_c^E$ the bounded component
which contains $E$ in its closure and abbreviate
$$\Sigma_c^E=\big\{(q,p) \in \Sigma_c: q \in \mathcal{K}_c^E\big\}.$$
If $\mu=0$ there is a single critical value of $U$ which can be
computed to be $-\frac{3}{2}$. For $c<-\frac{3}{2}$ the Hills region
has one bounded and one unbounded component. The bounded component
contains $E$ in its closure and we define $\mathcal{K}_c^E$ and
$\Sigma_c^E$ as before. For later reference let us abbreviate
$$\kappa=\left\{\begin{array}{cc}
U(L_1) & \textrm{if}\,\,\mu \in (0,1)\\
-\frac{3}{2} & \textrm{if}\,\,\mu=0
\end{array}\right.$$
to be the first critical value of $H$.

We denote by $\varphi\colon T^* \mathbb{R}^n \to T^*\mathbb{R}^n$
the symplectomorphism $(q,p) \mapsto (-p,q-E)$ and by $\iota \colon
T^* \mathbb{R}^n \to T^* S^n$ the inclusion obtained by interpreting
$\mathbb{R}^n$ as the chart of $S^n$ under stereographic projection.
We define the Moser regularization of $\Sigma^E_c$ to be
$$\overline{\Sigma}^E_c=\mathrm{cl}\big(\iota \varphi(\Sigma^E_c)\big)
\subset T^* S^n$$ where $\mathrm{cl}$ means closure. Recall that a
hypersurface $\Sigma \subset T^*S^n$ is called \emph{fiberwise
starshaped}, if for each $p \in S^n$ the intersection $\Sigma \cap
T_p^*S^n$ bounds a starshaped domain in $T^*_p S^n$. The following
result was proved in \cite{albers-frauenfelder-koert-paternain}.
\begin{thm}[Albers-Frauenfelder-van Koert-Paternain]\label{afkp}
In the planar case, i.e.\,$n=2$, if $c<\kappa$ then
$\overline{\Sigma}^E_c$ is a fiberwise starshaped hypersurface in
$T^*S^2$.
\end{thm}
\textbf{Remark: } Of course this holds also for the regularization
around the moon. Just replace $\mu$ by $1-\mu$ and the moon becomes
the earth while the earth becomes the moon.
\\ \\
\textbf{Remark: } We expect that the same result also holds in
higher dimensions, although we have not checked the details.

\section{The Birkhoff involution}\label{bi}

We write $q=(q_1,q_2,q^3)$ with $q^3 \in \mathbb{R}^{n-2}$ for a
vector in $\mathbb{R}^n$ The Hamiltonian $H$ of the circular
restricted three body problem is invariant under the antisymplectic
involution $B \colon T^* \mathbb{R}^n \to T^*\mathbb{R}^n$ given by
$$B(q_1,q_2,q^3,p_1,p_2,p^3)=(q_1,-q_2,-q^3,-p_1,p_2,p^3).$$
 Indeed, since $H$ is invariant under
$B$ but the symplectic form $\omega=\sum_{i=1}^n dq_i \wedge dp_i$
on $T^* N$ is antiinvariant, the Hamiltonian vector field $X_H$ of
$H$ is antiinvariant as well. In particular, if $w \in
C^\infty(\mathbb{R},T^* \mathbb{R}^n)$ is a trajectory of the
Hamiltonian flow of $H$, i.e.\,a solution of the ODE $\partial_t
w=X_H(w)$ and $Rw \in C^\infty(\mathbb{R},T^* \mathbb{R}^n)$ is
defined to by
$$Rw(t)=B(w(-t))$$
then $Rw$ is still a trajectory of the Hamiltonian flow. Periodic
orbits which are invariant under the involution $R$ are called
\emph{symmetric periodic orbits} and they played at least since the
work of Birkhoff \cite{birkhoff} a major role in the study of the
dynamics of the restricted three body problem.

We are now in position to proof the Observation we mentioned in the
Introduction.
\\ \\
\textbf{Proof of the Observation: } First we note that we can extend
$B$ to the Moser regularization as follows. Let $\rho \colon S^n \to
S^n$ be the reflection along the equator as in (\ref{refl}). The
adjoint of its differential $d^* \rho \colon T^* S^n \to T^* S^n$ is
a symplectic involution which commutes with the antisymplectic
involution $I \colon T^* S^n \to T^* S^n$ whose restriction to each
fibre is given by $I|_{T^*_p S^n}=-\mathrm{id}|_{T^*_p S^n}$. Their
composition $I \circ d^*\rho$ is an antisymplectic involution on
$T^* S^n$ which extends $B$. Thinking of $O(2)$ as the semidirect
product $O(2)=SO(2) \rtimes \mathbb{Z}_2$ we get an $O(2)$-action on
the free loop space of $T^* S^n$ as follows. For the subgroup
$SO(2)$ in $O(2)$ the action is just given by rotation of the
domain. If $r$ is the generator of $\mathbb{Z}_2$ and $w \in
\mathscr{L}_{T^* S^n}$, then the action $(\mathrm{id},r) \in SO(2)
\rtimes \mathbb{Z}_2$ is given by $(\mathrm{id},r)_* w=Rw$. If we
look at the shadow of $w$ in the free loop space of $S^n$ after
applying the footpoint projection $\pi \colon T^* S^n \to S^n$ we
precisely recover the twisted $O(2)$-action on $\mathscr{L}$ from
the introduction. This finishes the proof of the Observation. \hfill
$\square$
\\ \\
We finally prove the Corollary from the introduction.
\\ \\
\textbf{Proof of the Corollary: }By Theorem~\ref{afkp}
$\overline{\Sigma}_c^E$ bounds a Liouville domain $D_c$ in $T^*S^2$
which is isotopic to the unit disk bundle in $T^*S^2$. Hence by a
theorem of Abbondandolo-Schwarz, Salamon-Weber, and Viterbo
\cite{abbondandolo-schwarz,salamon-weber, viterbo} its symplectic
homology computes the loop space homology of $S^2$. Looking at the
$O(2)$-invariant symplectic homology \cite{bourgeois-oancea} we get
$$SH_*^{O(2)}(D_c)=H_*^{O(2)}(\mathscr{L}_{S^2})$$
where by the Observation the $O(2)$-action is the twisted
$O(2)$-action on the loop space of $S^2$. Lodder's computations
(\ref{lod}) now imply the Corollary. \hfill $\square$

\end{document}